
\input amstex

\mathsurround=2pt
\tolerance=10000
\def\mod{\text{mod}}
\def\fV{{\frak V}}
\def\vp{{\varphi}}

\def\]]{{\rbrack\!\rbrack}}
\def\[[{{\lbrack\!\lbrack}}

\def\ab0{{(a,b,0)}}

\def\fN{{\frak N}}

\def\BN{{\Bbb N}}
\def\BR{{\Bbb R}}

\def\cF{{\Cal F}}

\def\cC{{\Cal C}}
\def\cG{{\Cal G}}

\def\cN{{\Cal N}}

\def\cF{{\Cal F}}
\def\cD{{\Cal D}}

\def\cR{{\Cal R}}
\def\cP{{\Cal P}}
\def\cT{{\Cal T}}

\def\dim{\text{dim}}

\def\cC{{\Cal C}}
\def\cS{{\Cal S}}

\def\s1{{\textstyle {1\over 2}}}
\def\det{\text{det}}

\documentstyle{amsppt}
\NoBlackBoxes
\hoffset .5in
\topmatter

\title  COMPOSITE DIFFERENTIABLE FUNCTIONS\endtitle
\author Edward Bierstone, Pierre D. Milman and Wies\l aw Paw\l ucki\endauthor
\address Department of Mathematics, University of Toronto, Toronto, Canada M5S
1A1. \endaddress
\address Jagiellonian University, Institute of Mathematics, ul. Reymonta 4,
30-059 Krak\'ow, Poland. \endaddress

\thanks 1980 {\it Mathematics Subject Classification}
(1985 {\it Revision}).  Primary 32B20, 58C27; Secondary 32K15, 58C25.
\endthanks
\thanks{The first two authors' research was partially
supported by NSERC operating grants OGP 0009070 and OGP 0008949.}\endthanks
\endtopmatter

\document
\baselineskip=22pt
\head 1.  INTRODUCTION\endhead
\medskip

In the early 1940's, Whitney proved that every $\cC^\infty$
even function $f(x)$ can be written $f(x)=g(x^2)$,
where $g$ is $\cC^\infty$ [28].
About twenty years later, Glaeser (answering a question
posed by Thom in connection with the $\cC^\infty$ preparation
theorem) showed that a $\cC^\infty$ function $f(x)=f(x_1,\ldots,x_m)$
which is invariant under permutations of the coordinates
can be expressed $f(x)=g\big(\sigma_1(x),\ldots,\sigma_m(x)\big)$,
where $g$ is $\cC^\infty$ and the $\sigma_i(x)$ are the
elementary symmetric polynomials [10].
Of course, not every $\cC^\infty$ function $f(x)=f(x_1,\ldots,x_m)$
which is constant on the fibres of a (proper or semiproper) real
analytic mapping $y=\varphi(x)$, $y=(y_1,\ldots,y_n)$,
can be expressed as a composite $f=g\circ\varphi$, where
$g$ is $\cC^\infty$.
We will say that $\varphi$ has the {\it $\cC^\infty$ composite
function property} if every $\cC^\infty$ function $f(x)$
which is ``formally a composite with $\varphi$'' (see
Definition 1.1 below) can be written $f=g\circ\varphi$,
where $g(y)$ is $\cC^\infty$.
The theorem of Glaeser asserts that a semiproper real analytic
mapping $\varphi$ which is generically a submersion has the
$\cC^\infty$ composite function property.
The $\cC^\infty$ composite function property depends
only on the image $X$ of $\varphi$, which is a closed
subanalytic set [1] (cf. Corollary 1.5 below).
Bierstone and Milman have proved, more generally, that a
closed ``Nash subanalytic'' set $X$ has the $\cC^\infty$ composite
function property [1] (cf. [19,23,26]);
the class of Nash subanalytic sets includes all semianalytic
sets.
The $\cC^\infty$ composite function property is equivalent
to several other natural geometric and algebraic conditions
on a closed subanalytic set [6]; in particular, to a formal
semicoherence property (a stratified real version of the
coherence theory of Oka and Cartan).
Paw\l ucki has constructed an example of a closed subanalytic
set which is not semicoherent [20].
Thus the $\cC^\infty$ composite function property does not hold
in general, but distinguishes an important class of
subanalytic sets.

In this article, we introduce a new point of view towards
Glaeser's theorem, with respect to which we can formulate
a ``$\cC^k$ composite function property'' that is satisfied
by all semiproper real analytic mappings (Theorems 1.2 and 1.3
below).
As a consequence, we see that a closed subanalytic set $X$
satisfies the $\cC^\infty$ composite function property if and
only if the ring $\cC^\infty(X)$ of $\cC^\infty$ functions on $X$
is the intersection of all finite differentiability classes
(Corollary 1.5).

Let $k\in\BN\cup\{\infty\}$, where $\BN$ denotes the nonnegative
integers.
Suppose that $A$ is a locally closed subset of $\BR^n$ and that
$B\subset A$ is closed (in the relative topology of $A$).
We let $\cC^k(A;B)$ denote the Fr\'echet algebra of restrictions
to $A$ of $\cC^k$ real-valued functions that are defined on
neighbourhoods of $A$ and are $k$-flat on $B$.
(``$k$-flat on $B$'' means ``vanishing on $B$ together with all
partial derivatives of orders $\le k$''.)
$\cC^k(A)$ means $\cC^k(A;\emptyset)$.

Recall that a continuous mapping $\vp$: $M\to N$ is called
{\it semiproper} if, for each compact subset $K$ of $N$,
there is a compact subset $L$ of $M$ such that $\vp(L)=K\cap\vp(M)$.
Let $\vp$: $\Omega\to\BR^n$ denote a semiproper real analytic
mapping defined on an open subset $\Omega$ of $\BR^m$.
Then $X=\vp(\Omega)$ is a closed subanalytic subset of $\BR^n$.
(For the definitions and basic properties of subanalytic and Nash
subanalytic sets, see [1,5,7,8,9,14,16].)
Let $Z$ denote a closed subanalytic subset of $X$.
The mapping $\vp$ induces injective homomorphisms
$\vp^*$: $\cC^k(X;Z)\to \cC^k\big(\Omega;\vp^{-1}(Z)\big)$ given by
composition $\vp^*(g)=g\circ\vp$.

\medskip

\noindent {\bf Definition 1.1.}  {\it Let $\big(\vp^* \cC^k(X)\big)^\wedge$
denote the subalgebra of all functions $f\in \cC^k(\Omega)$ such that
$f$ is ``formally a composite with $\vp$''; i.e., for each $a\in X$,
there is $g\in \cC^k(X)$ such that $f-\vp^*(g)$ is $k$-flat
on $\vp^{-1}(a)$.
Set $\big(\vp^*\cC^k(X;Z)\big)^\wedge = \big(\vp^*\cC^k(X)\big)^\wedge
\cap \cC^k\big(\Omega;\vp^{-1}(Z)\big)$.}

\medskip

\proclaim{Theorem 1.2}  Assume that $X$ is compact.
Then, for each $k\in\BN$, there is an integer $\ell_k\ge k$
such that $\big(\vp^* \cC^\ell (X;Z)\big)^\wedge\subset
\vp^* \cC^k (X;Z)$ for all $\ell\ge\ell_k$.
\endproclaim
\medskip

Glaeser's theorem follows from Theorem 1.2:
Assume that $\vp$ is generically a submersion, so that the
interior $\text{int}\, X$ of $X$ is dense in $X$.
Suppose that $f\in\big(\vp^* \cC^\infty (X)\big)^\wedge$.
Let $g\in \cC^0(X)$ denote the unique function such that
$f=g\circ\vp$.
Then $g$ is $\cC^\infty$ in $\text{int}\, X$ and, by Theorem 1.2,
$g$ is the restriction to $X$ of a $\cC^k$ function, for all
$k\in\BN$.
Thus all partial derivatives of $g\mid\text{int}\, X$
extend continuously to $X$, and define a $\cC^\infty$ Whitney
field on $X$ (cf. \S3 below).
By Whitney's extension theorem, $g$ is the restriction
to $X$ of a $\cC^\infty$ function.

In general, put $\cC^{(\infty)}(X;Z)=\bigcap_{k\in\BN} \cC^k(X;Z)$
(and $\cC^{(\infty)}(X)=\cC^{(\infty)}(X;\emptyset)$).
The formula $\vp^*(g)=g\circ\vp$ defines an injective homomorphism
$\vp^*$:
$\cC^{(\infty)}(X;Z)\to \cC^\infty \big(\Omega;\vp^{-1}(Z)\big)$.
It is easy to see that $\bigcap_{\ell\in\BN}\big(
\vp^* \cC^\ell (X;Z)\big)^\wedge = \big(\vp^* \cC^\infty (X;Z)\big)^\wedge$,
so that Theorem 1.2 implies the following:

\medskip

\proclaim{Theorem 1.3}  $\big(\vp^* \cC^\infty (X;Z)\big)^\wedge =
\vp^* \cC^{(\infty)} (X;Z)$.
\endproclaim
\medskip

There is a continuous injection $\cC^\infty (X;Z)\hookrightarrow
\cC^{(\infty)}(X;Z)$, where the latter has a topology as the
inverse limit of the $\cC^k(X,Z)$, $k\in\BN$; it is easy to see
that $\cC^\infty (X;Z)$ is dense in $\cC^{(\infty)}(X;Z)$ (cf. [18, I.4.3]).
Moreover, $\big(\vp^* \cC^\infty (X;Z)\big)^\wedge$ is closed
in $\cC^\infty\big(\Omega;\vp^{-1}(Z)\big)$, so that Theorem 1.3
has the following corollaries:

\medskip

\proclaim{Corollary 1.4}  (cf. Tougeron [24]).
$\big(\vp^* \cC^\infty (X;Z)\big)^\wedge = \overline{\vp^* \cC^\infty
(X;Z)}$.
\endproclaim
\medskip

\proclaim{Corollary 1.5} The subalgebra $\vp^* \cC^\infty
(X;Z)$ is closed in $\cC^\infty \big(\Omega; \vp^{-1}(Z)\big)$ if and
only if $\cC^{(\infty)}(X;Z)=\cC^\infty (X;Z)$.
In particular, closedness of $\vp^* \cC^\infty (X)$ in $\cC^\infty (\Omega)$
depends only on the image $X$ of $\vp$ (cf. [1,6]).
\endproclaim
\medskip

\proclaim{Theorem 1.6} Let $X$ be a closed subanalytic
subset of $\BR^n$, and let $Z\subset X$ denote the points
which admit no neighbourhood in which $X$ is Nash.
(By [21], $Z$ is a closed subanalytic set.)
Then $\cC^{(\infty)}(X;Z) = \cC^\infty (X;Z)$.
\endproclaim

\medskip

Theorem 1.6 can be obtained by combining Corollary 1.5 with [6,
Theorem 3] and [22, Theorem 2].
In particular, if $X$ is a closed Nash subanalytic subset of
$\BR^n$, then $\cC^{(\infty)}(X) = \cC^\infty (X)$.
The same conclusion holds if $\text{dim}\, X\le 2$ or if
$X$ has pure codimension $1$ (by [22]).
It follows that, for any closed subanalytic subset $X$ of $\BR^4$,
$\cC^{(\infty)}(X) = \cC^\infty (X)$.
The example of [20] is a closed 3-dimensional subanalytic
subset $X$ of $\BR^5$ for which $\cC^{(\infty)}(X)\ne \cC^{\infty}(X)$.

\medskip

\noindent {\it Remark 1.7.}  Theorems 1.2, 1.3 and their corollaries
have more general module versions; i.e., versions that apply
not only to the solution of an equation of the form $f(x)=g\big(
\vp(x)\big)$ as above (where $g$ is the unknown), but also
more generally to systems of equations of the form $f(x)=A(x)\cdot
g\big(\vp(x)\big) + B(x)\cdot h(x)$, where $A(x)$, $B(x)$ are given
matrices of real analytic functions, $f(x)$ is a given
(vector-valued) $\cC^k$ function, and $g(y)$, $h(x)$ are the unknowns
(cf. [3,4]).

\medskip

\noindent {\it Remark 1.8.}  The results above can be stated, more
generally, for a semiproper mapping $\vp$: $M\to N$ of real analytic
manifolds:
The reduction to $N=\BR^n$ is immediate.
On the other hand, if $X=\vp (M)$ is compact, then there are
finitely many coordinate charts $\Omega_j$ in $M$ such
that $X=\cup\vp (\Omega_j)$ and $\vp$ is semiproper when regarded
as a mapping from the disjoint union $\Omega$ of the $\Omega_j$.
Therefore, we can reduce the statement
of Theorem 1.2 for $\vp$: $M\to N$ to the theorem as stated above.

\vskip .25in

\head 2.  DIVISION IN RINGS OF FORMAL POWER SERIES \endhead

\medskip

We will prove Theorem 1.2 by the method used in [3,4] based
on the division theorem of Grauert and Hironaka.

We totally order $\BN^n$ using the lexicographic ordering
of $(n+1)$-tuples
$(|\alpha|,\alpha_1,\break \ldots,\alpha_n)$, where
$\alpha=(\alpha_1,\ldots,\alpha_n)\in\BN^n$ and
$|\alpha|=\alpha_1+\cdots+\alpha_n$.
For every formal power series $F=\sum_\alpha f_\alpha y^\alpha\in
\BR[[y]] = \BR[[y_1,\ldots,y_n]]$, we define the {\it support}
$\text{supp}\, F = \{\alpha\in\BN^n:\, f_\alpha\ne 0\}$,
the {\it initial exponent} $\text{exp}\, F =\text{min}\,\text{supp}\,
F$ ($\text{ exp}\, F=\infty$ if $F=0$), and the {\it initial
monomial} $\text{ mon}\, F = f_\alpha y^\alpha$, where
$\alpha=\text{ exp}\, F$ (when $F\ne 0$).

\medskip

\proclaim{Theorem 2.1} (cf. [11,15]).
Let $F_1,\ldots,F_s\in\BR[[y]]\backslash\{0\}$, and let
$\alpha^i=\text{ exp}\, F_i$ $(i=1,\ldots,s)$.
Put $\Delta_0=\emptyset$, $\Delta_i=(\alpha^i+\BN^n)\backslash
\big(\bigcup_{j<i}\Delta_j\big)$ for $i=1,\ldots,s$, and
$\Delta=\BN^n\backslash\big(\bigcup_{i=0}^s\Delta_i\big)$.
Let $G\in\BR[[y]]$.
Then there are unique $Q_i$, $R\in\BR[[y]]$ $(i=1,\ldots,s)$
such that $\alpha^i+\text{ supp}\, Q_i\subset\Delta_i$ $(i=1,\ldots,s)$,
$\text{ supp}\, R\subset\Delta$, and
$$
G\ =\ \sum_{i=1}^s Q_i F_i +R \ .
$$
Moreover, $\text{ exp}\, R\ge \text{ exp}\, G$ and, for each $i$,
$\alpha^i + \text{ exp}\, Q_i\ge \text{ exp}\, G$.
\endproclaim

See [4, p. 207] for the proof.

\medskip

Let $I$ be an ideal in $\BR[[y]]$.
The {\it diagram of initial exponents} of $I$ is defined
as $\fN(I)=\big\{\text{ exp}\, F:\, F\in I\backslash\{0\}\big\}$.
Clearly, $\fN(I)+\BN^n = \fN(I)$.
This property implies that there is a smallest finite set
$\fV(I)\subset \fN(I)$ such that $\fV(I)+\BN^n=\fN(I)$.
The elements of $\fV(I)$ are called the {\it vertices} of
$\fN(I)$.

Let $\fV(I)=\{\alpha^1,\ldots,\alpha^s\}$, and let
$\{\Delta_i,\Delta\}$ denote the decomposition of $\BN^n$
determined by the $\alpha^i$, as in Theorem 2.1.

\medskip

\proclaim{Proposition 2.2} [4, Corollary 6.8].
Choose $F_1,\ldots,F_s\in I$ such that $\text{ exp}\, F_i =
\alpha^i$ $(i=1,\ldots,s)$.
Then $F_1,\ldots,F_s$ generate $I$.
Moreover there is a unique system $\{ G_1,\ldots, G_s\}$
of generators of $I$ such that, for each $i$, $\text{ supp}
(G_i-y^{\alpha^i})\subset\Delta$;
in particular, $\text{ mon}\, G_i=y^{\alpha^i}$.
\endproclaim
\medskip

We call $\{ G_1,\ldots, G_s\}$ the {\it standard basis} of $I$.
Using Theorem 2.1, we get:

\medskip

\proclaim{Proposition 2.3} (cf. [4, Corollary 6.9]).
Let $(y)$ denote the maximal ideal of $\BR[[y]]$.
Let $p\in\BN$.
Then the set $\{ y^\beta:\, \beta\in\Delta,|\beta|\le p\}$
is a basis of an $\BR$-linear complement of $I+(y)^{p+1}$
in $\BR[[y]]$.
In particular, if $(y)^{p+1}\subset I$, then $\{ y^\beta:\, \beta\in
\Delta\}$ is a basis of an $\BR$-linear complement of $I$
in $\BR[[y]]$.
\endproclaim
\medskip

Let $\cD(n)=\{\fN\subset\BN^n:\, \fN+\BN^n=\fN\}$.
We totally order $\cD(n)$ as follows:
To each $\fN\in\cD(n)$, associate the sequence $v(\fN)$ obtained
by listing the vertices of $\fN$ in ascending order and completing
the list to an infinite sequence by using $\infty$ for all the
remaining terms.
If $\fN_1,\fN_2\in\cD(n)$, we say that $\fN_1 < \fN_2$ provided
that $v(\fN_1)< v(\fN_2)$ with respect to the lexicographic
ordering on the set of all such sequences.

\vskip .25in

\head 3.  WHITNEY FIELDS\endhead

\medskip

Let $p\in\BN \backslash\{0\}$ and let $A$ be a locally closed subset of
$\BR^n$.
A $\cC^p$-{\it Whitney field} on $A$ is a polynomial
$F(a,y)=\sum_{|\alpha|\le p} {1\over\alpha!} F_\alpha (a) y^\alpha
\in \cC^0(A)[y]=\cC^0(A)[y_1,\ldots,y_n]$
which fulfills the following conditions:
$$
\left|\left(\partial^{|\beta|}F/\partial y^\beta\right)(a,0)-
\big(\partial^{|\beta|} F/\partial y^\beta\big)(b,a-b)\right| \big/
|a-b|^{p-|\beta|}\ \longrightarrow \ 0\ ,
$$
when $a\to c$, $b\to c$, $a\ne b$, for each $c\in A$ and $\beta\in\BN^n$
such that $|\beta|\le p$.

For any polynomial
$F(a,y)=\sum_{|\alpha|\le p} {1\over\alpha!} F_\alpha(a) y^\alpha\in
\cC^0(A)[y]$,
we set $F^k(a,y)=\sum_{|\alpha|\le k} {1\over\alpha!} F_\alpha(a)y^\alpha$,
where $k\le p$, $k\in\BN$.

For any $\cC^p$-mapping $\vp=(\vp_1,\ldots,\vp_n)$: $D\to\BR^n$ defined on an
open
subset $D$ of $\BR^m$, we denote by $T^p\vp$ its {\it Taylor field}
$$
T^p\vp(b,x) = \sum_{|\alpha|\le p} {1\over\alpha!} D^\alpha \vp(b) x^\alpha
\ = \ \left(\sum_{|\alpha|\le p} {1\over\alpha!} D^\alpha \vp_1(b)
x^\alpha,\ldots,
\sum_{|\alpha|\le p} {1\over\alpha!} D^\alpha \vp_n(b) x^\alpha\right)
\ .
$$
Put $\widetilde T^p\vp(b,x)=T^p\vp(b,x)-\vp(b)$
(the field of Taylor polynomials without constant terms).

\medskip

\proclaim{Lemma 3.1}
Let $\Lambda$ be a $\cC^1$-submanifold of $\BR^n$.
Let
$F(a,y)=\sum_{|\alpha|\le p} {1\over\alpha!} F_\alpha(a)y^\alpha\in
\cC^0(\Lambda)[y]$,
where
$F^{p-1}(a,y)\in \cC^1(\Lambda)[y]$.
Let $\vp:\, D\to\BR^n$ be a $\cC^p$-mapping
as above, and let $\Gamma$ be a $\cC^1$-submanifold of $D$ such that
$\vp(\Gamma)\subset\Lambda$.
Define a polynomial
$G\in \cC^0(\Gamma)[x]=\cC^0(\Gamma)[x_1,\ldots, x_m]$ of degree $\le p$, by
the
formula
$$
G(b,x)\ =\ F\big(\vp(b),\widetilde T^p\vp(b,x)\big)\ \mod\  (x)^{p+1}\  .
$$
Then $G^{p-1}\in \cC^1(\Gamma)[x]$ and,
for each $b\in\Gamma$ and $v\in T_b\Gamma$, we have
$$
\eqalign{
&D_{b,v} G^{p-1}(b,x)-D_{x,v} G(b,x)\cr
&=\
D_{a,u} F^{p-1}\big(\vp(b),\widetilde T^p\vp(b,x)\big)-
D_{y,u}F\big(\vp(b),\widetilde T^p\vp(b,x)\big)\ \mod\  (x)^p ,\cr}
$$
where $u=d_b\vp(v)$, and $D_{b,v}$ and $D_{x,v}$ stand for the
(directional) derivatives with respect to $b$ and $x$,
respectively, along the tangent vector $v$.
\endproclaim

This is a straightforward calculation.

\medskip

\proclaim{Proposition 3.2} (cf. [4, \S10]).
Let $\Lambda$ be a $\cC^p$-submanifold of $\BR^n$ and let
$F\in \cC^0(\Lambda)[y]$,
$\text{ deg}\, F\le p$. Then $F$ is a $\cC^p$-Whitney
field on $\Lambda$ if and only if $F^{p-1}\in \cC^1(\Lambda)[y]$ and
$D_{a,u} F^{p-1}(a,y)=D_{y,u}F(a,y)$, for every $a\in\Lambda$ and
$u\in T_a\Lambda$.
\endproclaim

\demo{Proof}
It is enough to consider the case $\Lambda=\BR^k\times 0$;
the proposition then follows easily.  $\square$
\enddemo

\medskip

\noindent {\it Remark 3.3.} [25, IV.2.5].
Let $U$ denote an open subset of $\BR^n$.
Let $u,v\in U$, and suppose that $\sigma$ is a
rectifiable curve in $U$ joining $u,v$.
Let $g\in \cC^k (U)$, $k\ge 1$.
It follows from the Mean Value Theorem (using an
approximation of $\sigma$ by piecewise-linear curves)
that
$$
|g(u)-g(v)|\ \le\ \sqrt{n} |\sigma| \sup_{a\in\sigma\atop
1\le j\le n} \left| {\partial g\over\partial y_j} (a)\right|\ ,
$$
where $|\sigma|$ denotes the length of $\sigma$.
If $g$ is $(k-1)$-flat at $v$, then (by iterating the
inequality above) we get
$$
|g(u)|\ \le\ n^{k/2} |\sigma|^k \sup_{a\in\sigma\atop
|\alpha|=k} \left| {\partial^{|\alpha|} g\over
\partial y^\alpha} (a)\right|\ .
$$
Now suppose that $F$ is a $\cC^\ell$ Whitney field on $\sigma$.
Let $\alpha\in\BN$, $|\alpha|<\ell$, and apply the preceding
inequality with $k=\ell-|\alpha|$ and
$\displaystyle{g(y)={\partial^{|\alpha|}\over\partial y^\alpha}} \big( f(y)-
F(v,y-v)\big)$, where $f\in \cC^\ell(U)$ is a Whitney
extension of $F$.
We obtain
$$
\left| {\partial^{|\alpha|}F\over \partial y^\alpha}
(u,0) - {\partial^{|\alpha}| F\over\partial y^\alpha}
(v,u-v)\right|\ \le\ n^{{\ell-|\alpha|\over 2}}
|\sigma|^{\ell-|\alpha|} \sup_{a\in\sigma\atop|\beta|=\ell}
|F_\beta(a)-F_\beta(v)|\ .
$$
(This inequality holds trivially if $|\alpha|=\ell$).

\medskip

Let $r\in\BN\backslash\{0\}$.
A compact subset $A$ of $\BR^n$ is called $r$-{\it regular}
if there is a constant $C>0$ such that any two points
$u,v\in A$ can be joined by a rectifiable curve $\sigma$
in $A$ of length $|\sigma|\le C|u-v|^{1/r}$.
The following is a version of l'H\^opital's rule
or Hestenes's lemma for $r$-regular sets (generalizing [27]).

\medskip

\proclaim{Proposition 3.4}  Let $A\supset B$
be compact subsets of $\BR^n$, where $A$ is $r$-regular.
Let $k\in\BN$ and let $F(a,y)=\sum_{|\alpha|\le kr}
{1\over\alpha!} F_\alpha(a)y^\alpha \in \cC^0
(A)[y]$, $y=(y_1,\ldots,y_n)$.
Assume that $F$ restricts to $\cC^{kr}$ Whitney fields
on $B$ and on $A\backslash B$.
Then $F^k(a,y)=\sum_{|\alpha|\le k}{1\over\alpha!}
F_\alpha (a) y^\alpha$ is a $\cC^k$ Whitney field on $A$.
\endproclaim

\demo {Proof}  For all $u,v\in A$ and each $\alpha\in\BN^n$,
$|\alpha|\le k$,
$$
\align
& {\partial^{|\alpha|} F^k\over \partial y^\alpha} (u,0) -
{\partial^{|\alpha|} F^k\over \partial y^\alpha} (v,u-v) \\
&=\ {\partial^{|\alpha|} F\over \partial y^\alpha}
(u,0) - {\partial^{|\alpha|} F\over\partial y^\alpha}
(v,u-v) + \sum_{k-|\alpha|<|\beta|\le kr-|\alpha|}
{1\over\beta!} F_{\alpha+\beta} (v) (u-v)^\beta\ .\tag 3.5
\endalign
$$
We can assume that $F$ is $kr$-flat on $B$ (by
Whitney's extension theorem).
Let $u,v\in A$, and let $\sigma$ be a rectifiable curve in
$A$ of length $\le C|u-v|^{1/r}$ joining $u,v$.
Write
$$
\mu_F (\sigma)\ =\ \sup_{a,b\in\sigma\atop |\beta|=kr}
|F_\beta (a) - F_\beta (b)|\ ,
$$
so that $\mu_F(\sigma)\to 0$ as $|u-v|\to 0$.
\enddemo

Case 1.  $\sigma\cap B =\emptyset$.
By Remark 3.3, for all $|\alpha|\le kr$,
$$
\eqalign{
\left| {\partial^{|\alpha|}F\over\partial y^\alpha} (u,0) -
{\partial^{|\alpha|} F\over\partial y^\alpha} (v,u-v)\right|\
&\le \ n^{{kr-|\alpha|\over 2}} |\sigma|^{kr-|\alpha|}
\mu_F (\sigma)\cr
&\le \ c_1 |u-v|^{k-|\alpha|} \mu_F(\sigma)\ ,\cr}
$$
where $c_1$ is independent of $u,v$.
By (3.5), if $|\alpha|\le k$, then
$$
\left| {\partial^{|\alpha|}F^k\over\partial y^\alpha} (u,0) -
{\partial^{|\alpha|} F^k\over\partial y^\alpha} (v,u-v)\right|
\ \le \ c_2 |u-v|^{k-|\alpha|} (\mu_F(\sigma)+|u-v|)\ .
$$

Case 2.  $\sigma\cap B\ne\emptyset$.
If $w\in\sigma$, let $\sigma_{u,w}$ denote the shortest
part of $\sigma$ joining $u,w$.
Then there exists $w\in\sigma$ such that $\sigma_{u,w}\cap B =
\{ w\}$.
Using the Mean Value Theorem (as in Remark 3.3), we see that
if $|\alpha|\le kr$, then
$$
\eqalign{
|F_\alpha(u)|\ &\le\ n^{{kr-|\alpha|\over 2}} |\sigma_{u,w}|^
{kr-|\alpha|} \sup_{a\in\sigma_{u,w}\atop |\beta|=kr}
|F_\beta(a)|\cr
&\le\ n^{{kr-|\alpha|\over 2}} |\sigma|^{kr-|\alpha|}
\mu_F (\sigma)\cr
&\le\ c_1 |u-v|^{k-|\alpha|} \mu_F(\sigma)\ ,\cr}
\leqno (3.6)
$$
and (applying (3.6) with $v$ in place of $u$)
$$
\left| \sum_{|\beta|\le kr-|\alpha|} {F_{\alpha+\beta}(v)\over
\beta!} (u-v)^\beta\right|\ \le\ c_3 |u-v|^{k-|\alpha|}
\mu_F (\sigma)\ .
$$
Thus
$$
\left| {\partial^{|\alpha|} F\over \partial y^\alpha}
(u,0) - {\partial^{|\alpha|} F\over \partial y^\alpha}
(v,u-v)\right|\ \le\ (c_1 + c_3) |u-v|^{k-|\alpha|}
\mu_F(\sigma)\ ,
$$
and the required estimate again follows from (3.5).  $\square$

\vskip .25in

\head 4. STRATIFICATION OF A SUBANALYTIC MAPPING\endhead

\medskip

We will use a theorem of Hardt on stratification
of mappings [12,13] in the version of \L ojasiewicz
[17].

A {\it subanalytic leaf} in $\BR^m$ means a connected
subanalytic subset of $\BR^m$ which is an analytic submanifold.
Let $E$ be a subanalytic subset of $\BR^m$. A {\it subanalytic
stratification} of $E$ (in $\BR^m$) is a partition $\cS$
of $E$ into subanalytic leaves in $\BR^m$, called {\it
strata}, locally finite in $\BR^m$, such that, for each
$S\in \cS$, the boundary $(\overline S\backslash S)\cap E$ is a union
of strata of dimension $< \dim S$. A mapping
$\vp$: $E\to \BR^n$ is called {\it subanalytic} if its graph is
subanalytic in $\BR^m\times\BR^n$. A partition $\cP$ of
$E$ is called {\it compatible} with a family $\cF$ of
subsets of $E$ if, for each $P\in\cP$ and $F\in\cF$, either
$P\subset F$ or $P\subset E\backslash F$.

\medskip

\proclaim{Theorem 4.1} (cf. [17]).
Let $\vp$: $E\to\BR^n$ be a
continuous subanalytic mapping defined on a compact subanalytic
subset $E$ of $\BR^m$. Put $X=\vp(E)$. Let $\cF$ and $\cG$
be finite families of subsets of $E$ and $X$ which are subanalytic in
$\BR^m$ and $\BR^n$, respectively. Then there exist finite
subanalytic stratifications $\cS$ and $\cT$ of $E$ and
$X$, respectively, such that:

(1)  For each $S\in \cS$, $\vp(S)\in\cT$ and there
is a commutative diagram
$$

\matrix
S&&{\smash{\mathop{\longrightarrow}\limits^{h}}}&T\times P\\
\ {\scriptstyle{\vp|S}}&\searrow&&\!\!\!\!{\swarrow\rlap{
$\vcenter{\hbox{$\scriptstyle\pi$}}$}}\\
&&T\endmatrix
$$
where $T=\vp(S)$, $P$ is a bounded subanalytic leaf in
$\BR^s$ for some $s$, $h$ is an analytic subanalytic isomorphism and
$\pi$ denotes the natural projection.

(2)  $\cS$ is compatible with $\cF$ and $\cT$ is
compatible with $\cG$.
\endproclaim

\medskip

A pair $(\cS, \cT)$ of subanalytic stratifications as in
(1) will be called a {\it stratification} of $\vp$. If we
weaken this definition by allowing $\cS$ to be any finite
partition into subanalytic leaves (but still requiring $\cT$ to
be a stratification), then the pair $(\cS,\cT)$ will be
called a {\it semistratification}\/ of $\vp$.

\medskip

\noindent {\it Remark 4.2.}
Clearly, Theorem 4.1 is true for a subanalytic mapping
$\vp$: $E\to\BR^n$
which is defined  on a bounded subanalytic subset $E$ of $\BR^m$ and extends
continuously to $\overline E$.

\medskip

\noindent {\it Remark 4.3.}
If $(\cS,\cT)$ is  a semistratification of
$\vp$ and $Y$ is a subanalytic subset of $X$ such that
$\cT'=\{T\cap Y:\, T\in \cT\}$ is a stratification of $Y$, then the
pair $(\cS',\cT')$, where
$\cS'=\{S\cap\vp^{-1}(Y):\, S\in\cS\}$, is a
semistratification of $\vp|\vp^{-1}(Y):\,\vp^{-1}(Y)\to Y$.

\medskip

\noindent {\it Remark 4.4.}
If $(\cS,\cT)$ is a semistratification of $\vp$,
$S\in\cS$, $T\in\cT$ and $\vp(S)=T$, then there exists a
subanalytic leaf $\Gamma\subset S$ such that $\vp|\Gamma$: $\Gamma\to T$
is an analytic isomorphism.

\medskip

Let $\vp$: $E\to \BR^n$ be a bounded subanalytic mapping defined on a bounded
subanalytic subset $E$ of $\BR^m$.
For each $q\in\BN$, $q\ge 1$, we define the $q$-{\it fold fibre
product of} $E$ {\it with respect to} $\vp$,
$$
E_\vp^q=\left\{ \underline b=(b_1,\ldots, b_q)\in E^q:\
\vp(b_1)=\cdots=\vp(b_q) \right\}\ ;
$$
$E_\vp^q$ is a subanalytic subset of $(\BR^m)^q$.
There is a
natural mapping $\Phi$: $E_\vp^q\to\BR^n$ defined by
$\Phi(b_1,\ldots,b_q)=\vp(b_1)$.
Suppose that $(\cS,\cT)$ is a
semistratification of $\vp$. Let $\cS^{(q)}$ denote the family of all
non-empty sets of the form $(S_1\times\cdots\times S_q)\cap E_\vp^q$,
where $S_1,\ldots S_q \in\cS$. It is easy to see (using
Remark 4.3) that
$(\cS^{(q)},\cT)$ is a semistratification of $\Phi$.

\vskip .25in

\head 5.  IDEALS OF RELATIONS \endhead

\medskip

Let $E$ be a bounded open subanalytic subset of
$\BR^m$, and let $\vp =(\vp_1,\ldots,\vp_n): \ E\to\BR^n$ be
a mapping which extends to be analytic
in a neighbourhood of $\overline E$. Let
$p\in\BN$, $p\ge 1$.
Put $X=\vp(E)$ and $q={n+p\choose p}$.

For each $a\in X$, we define the {\it ideal of relations of order $p$
among $\vp_1,\ldots,\vp_n$ over} $a$ as
$$
\cR_p(a)=\{W\in\BR[[y]]:\
W\big(\widetilde T^p\vp(b,x)\big)\ =\ 0\ \mod\ (x)^{p+1},\
\hbox{for each}\ b\in\vp^{-1}(a)\}
$$
Clearly, $(y)^{p+1}\subset \cR_p(a)$.
Put
$$
L_\beta^\alpha(b) = D_x^\beta\big((\widetilde T^p\vp(b,x))^\alpha\big)(0)\  ,
$$
where $\alpha\in\BN^n$, $\beta\in\BN^m$, $|\alpha|\le p$,
$|\beta|\le p$, $b\in E$.
($D_x^\beta$ denotes the partial derivative $\partial^{|\beta|}/
\partial x^\beta$ with respect to the variables $x$.)
Then for every formal power series
$W=\sum_{\alpha} W_\alpha y^\alpha\in\BR[[y]]$,
$W\in \cR_p(a)$ if and only if $\{W_\alpha\}_{|\alpha|\le p}$
satisfies the following system of linear equations:
$$
\sum_{|\alpha|\le p} W_\alpha \cdot L_\beta^\alpha(b)\ =\ 0
\qquad \big(|\beta|\le p,\ b\in\vp^{-1}(a)\big) \ .
\leqno(5.1)
$$
Let $(\cS,\cT)$ be a semistratification of $\vp$. This induces a
semistratification $(\cS^{(q)},\cT)$ of the mapping
$\Phi$: $E_\vp^q\to X$.
Let $l\in\BN$, $1\le l\le p$. Denote by
$\pi_l$: $\BR[[y]]\to\BR[y]$ the projection
$\pi_l\left(\sum_\alpha W_\alpha y^\alpha\right) = \sum_{|\alpha|\le l}
W_\alpha
y^\alpha$.
It is clear that $\pi_l(\cR_r(a))\subset \pi_l(\cR_p(a))\subset \c
R_l(a)$, when $l\le p\le r$ and $a\in X$.

For each $\underline{b}=(b_1,\ldots,b_q)\in E_\vp^q$, let
$\rho^0(\underline{b})$
denote the rank of the matrix of the system
$$
\sum_{|\alpha|\le p} W_\alpha L_\beta^\alpha(b_\nu) =0 \qquad
(|\beta|\le p,\ \nu=1,\ldots,q)\ ,
\leqno(5.2)
$$
and $\rho^1(\underline{b})$ the rank of the matrix of the system
$$
\sum_{l<|\alpha|\le p} W_\alpha L_\beta^\alpha(b_\nu) \ =\ 0\qquad
(|\beta|\le p,\ \nu=1,\ldots,q)\ .
\leqno(5.3)
$$
(Of course, $\rho^0(\underline b )$ depends on $p$,
and $\rho^1(\underline b )$ on $p$ and $l$.)
For each $T\in\cT$, we put
$$
\eqalign{
&\sigma_T^0\ =\ \max\{\rho^0(\underline{b}): \
\underline{b}\in\Phi^{-1}(T)\}\ ,\cr
&\sigma_T^1\ =\ \max\{\rho^1(\underline{b}): \
\underline{b}\in\Phi^{-1}(T)\}\ .\cr}
$$
Observe that if $a\in T$ is a point such that
$$
\eqalign{
&\max\{\rho^0(\underline{b}):\
\underline{b}\in\Phi^{-1}(a)\}\ =\ \sigma_T^0\ ,\cr
&\max\{\rho^1(\underline{b}): \
\underline{b}\in\Phi^{-1}(a)\}\ =\ \sigma_T^1\ ,\cr}
$$
then
$\dim\pi_p(\cR_p(a))=q-\sigma_T^0$ and
$\dim\pi_l(\cR_p(a))=q-\sigma_T^0-(q-{l+n\choose n}-\sigma_T^1)=
{l+n\choose n}+\sigma_T^1-\sigma_T^0$.

For each $T\in\cT$, let us take leaves $S_T^0, S^1_T\in\cS^{(q)}$
such that $\Phi(S_T^0)=\Phi(S_T^1)=T$ and
$$
\eqalign{
&\max\{\rho^0(\underline{b}):\ \underline{b}\in S_T^0\}\ =\ \sigma_T^0
\ ,\cr
&\max\{\rho^1(\underline{b}):\ \underline{b}\in S_T^1\}\ =\ \sigma_T^1
\ .\cr}
$$
Then
there is a closed nowhere dense subset $Z_T^{pl}$ of
$T$, subanalytic in $\BR^n$, such that, for each $a\in T\backslash
Z_{T}^{pl}$,
$$
\eqalign{
&\max\{\rho^0(\underline{b}): \ \underline{b}\in S_T^0\cap\Phi^{-1}(a)\}\ =\
\sigma_T^0\ ,\cr
&\max\{\rho^1(\underline{b}):\ \underline{b}\in
S_T^1\cap\Phi^{-1}(a)\}\ =\ \sigma_T^1\ .\cr}
$$
Put $\omega_T^{pl}={l+n\choose n}+\sigma_T^1-\sigma_T^0$. Then
$\omega_T^{pl}=\dim \pi_l(\cR_p(a))$, for $a\in T\backslash Z_T^{pl}$.
We have $\omega_T^{rl}\le \omega_T^{pl}$, when $l\le p\le r$.
Therefore, for
every $l$, there exists $p_l\ge l$ such that $\omega_T^{pl}$ is constant
for $p\ge p_l$. Since $\cT$ is finite, we can make $p_l$ independent of
$T$. This gives the following:

\medskip

\proclaim{Lemma 5.4} For every positive $l\in\BN$,
there exists $p_l\in \BN$,
$p_l\ge l$, such that if  $p>p_l$ and $T\in \cT$, then
$\pi_l(\cR_p(a))=\pi_l(\cR_{p-1}(a))$, for each
$a\in T\backslash\left(Z_T^{pl}\cup Z_T^{(p-1)l}\right)$.
\endproclaim

\medskip

Now let us fix any $p>p_l$. Take $T\in\cT$.
The diagram of initial exponents $\fN(\cR_p(a))$, $a\in
T\backslash Z_T^{pl}$, takes only finitely many values;
choose $a\in T\backslash Z_T^{pl}$ with the minimum
value (cf. \S2).
Set
$\Delta(a)=\BN^n\backslash\cN(\cR_p(a))$. There exists
$\underline{b}\in S_T^0\cap\Phi^{-1}(a)$ such that
$\rho^0(\underline{b})=\sigma_T^0$.
By Proposition 2.3, for any such $\underline{b}$, we have
$\#\Delta(a)=\rho^0(\underline{b})=\sigma^0_T=\text{ rank}
\{L_\beta^\alpha(b_\nu)\}$ (where
$\alpha\in\Delta(a)$, $|\beta|\le p$, $\nu=1,\ldots,q$).
There exists a nonzero determinant
$$
M_T(\underline{b})=\det\left(\{L_{\beta_k}^\alpha(b_{\nu_k})\}
\right)
\quad(\alpha\in\Delta(a),\ k=1,\ldots,\sigma_T^0)\ .
$$
Put $\Theta_T^{pl} = \{\underline{c}\in S^0_T:\  M_T(\underline{c})=0\}$ and
$\Sigma_T^{pl}=T\backslash\Phi\left(S_T^0\backslash\Theta_T^{pl}\right)$.
Clearly, $\Theta_T^{pl}$ and $\Sigma_T^{pl}$ are nowhere dense in $S_T^0$
and in $T$, and are subanalytic in $(\BR^m)^q$ and in $\BR^n$, respectively.
Since $M_T(\underline{b})\ne 0$, it follows from Cramer's rule that the
system (5.2), and therefore the system (5.1), is equivalent
to the system
$$
W_\gamma - \sum_{\alpha\not\in\Delta(a)} W_\alpha\cdot
N_{T,\gamma}^\alpha(\underline{b})/M_T(\underline{b})\ =\
0\quad(\gamma\in\Delta(a)),
\leqno(5.5)
$$
where $N_{T,\gamma}^\alpha(\underline{b})$ are appropriate minors of the
matrix of the system (5.2).
Let $\{\alpha^1,\ldots\alpha^s\}=\fV(\cR_p(a))$, where
$\alpha^1<\ldots<\alpha^s$, and let
$r=\max\{i:\ |\alpha^i|\le p\}$. Since $(y)^{p+1}\subset\cR_p(a)$,
$|\alpha^i|=p+1$ for $i>r$.
Observe that the standard basis of $\cR_p(a)$ consists of $y^{\alpha^i}$
($i=r+1,\ldots,s$) and
$$
G_i(y)=y^{\alpha^i}+\sum_{\gamma\in\Delta(a)}
\left(N_{T,\gamma}^{\alpha^i}(\underline{b}))
/M_T(\underline{b})\right)y^\gamma \quad (i=1,\ldots,r).
$$
Moreover, we have
$$
y^{\alpha^i}
+\sum_{\gamma\in\Delta(a)}\left(N_{T,\gamma}^{\alpha^i}(\underline{b'}))
/M_T(\underline{b'})\right)y^\gamma \in\cR_p(\Phi(\underline{b'})),
$$
for any $\underline{b'}\in S_T^0\backslash\Theta_T^{pl}$
such that $\Phi(\underline{b'})\not\in Z_T^{pl}$, and $i=1,\ldots,r$.
Suppose that $a'\in T\backslash (Z_T^{pl}\cup \Sigma_T^{pl})$.
Since $\fN(\cR_p(a))\le\fN(\cR_p(a'))$, it follows that
$\alpha^i=\exp G_i\in\fN(\cR_p(a'))$ for each $i$,
so that $\fN(\cR_p(a))\subset\fN(\cR_p(a'))$ and hence
they are equal.
In other words,
$\Delta(a')$ is
independent of $a'\in T\backslash\left(Z_T^{pl}\cup\Sigma_T^{pl}\right)$. Let
us
denote it by $\Delta_T$.

Now taking an appropriate semistratification of $\Phi$ and using
Remark 4.4,
we can find a finite subanalytic stratification $\{\Lambda_j\}$ ($j\in J$) of
$X$, compatible with the sets $T$, $Z_T^{pl}$, $Z_T^{(p-1)l}$,
$\Sigma^{pl}_T$, and a family $\{\Gamma_j\}$ ($j\in J$) of subanalytic leaves
in $E_\vp^q$ such that, for each $j\in J$,
$\Phi|\Gamma_j:\ \Gamma_j\to\Lambda_j$ is an analytic isomorphism, and
$\Gamma_j\subset S_T^0\backslash \Theta_T^{pl}$, whenever $\Lambda_j\subset
T\backslash \Sigma_T^{pl}$.

Suppose that $j\in J$, $T\in \cT$, $\Lambda_j\subset T$ and
$\dim\Lambda_j=\dim T$. Then $\Lambda_j$ is open in $T$ and
$\Lambda_j\subset T\backslash\left(Z_T^{pl}\cup\Sigma_T^{pl}\right)$. We put
$M_j=M_T$ and $\Delta_j=\Delta_T$ in this case.

 \vskip .25in
\head 6. PROOF OF THEOREM 1.2 \endhead
\medskip

Since $\vp$ is semiproper and $X=\vp(\Omega)$ is compact,
there exists a bounded open subanalytic subset $E\subset \BR^m$
such that $\overline E\subset\Omega$, $\vp(E)=X$ and $\vp|E:\ E\to\BR^n$
is semiproper.
Theorem 1.2 reduces to the following proposition (by induction
on dimension).

\medskip

\proclaim{Proposition 6.1}
Let $A$ be a closed subset of $X$, subanalytic in $\BR^n$,
of dimension $d$.
Let $k\in\BN$. Then there exists a closed subset
$A'$ of $A$ and an integer $t\ge k$ such that:

\text{ (1)} $A'$ is subanalytic of dimension $<d$.

\text{ (2)} If $f\in\cC^t(\Omega)$ is $t$-flat on $\vp^{-1}(A'\cup Z)$ and, for
each
$a\in A$, there exists a polynomial $W_a\in\BR[y]$ such that
$T^tf(b,x)=W_a(\widetilde T^t\vp(b,x))\ \mod\  (x)^{t+1}$, for all
$b\in\vp^{-1}(a)\cap E$, then there exists $g\in\cC^k(\BR^n;Z)$ such
that $f-g\circ \vp$ is $k$-flat on $\vp^{-1}(A)\cap E$.
\endproclaim

\demo{Proof}  Let $r$ denote a positive integer such that
every connected component of $A\cup Z$ is $r$-regular
[5, Th. 6.10].
Take any integer $l\ge kr$.
Let $(\cS,\cT)$ be a (semi) stratification of $\vp|E$ such
that $\cT$ is compatible with $A$ and $Z$.
We will use the results and notation of \S5.
Fix an integer $p>p_l$ and
take the stratification $\{\Lambda_j\}$ ($j\in J\}$ of $X$ and the family
$\{\Gamma_j\}$ ($j\in J$) of subanalytic leaves in $E^q_\vp$ as at the end
of Section~5. Observe that the mapping $\Phi$: $E_\vp^q\to X$ extends to
$\Phi$: $\Omega^q\to X$, and $\overline{E_\vp^q}\subset\Omega^q$.
For each $j\in J$, $\Phi\left(\overline\Gamma_j\backslash
\Gamma_j\right)=\overline\Lambda_j\backslash\Lambda_j$, because
$\Phi|\Gamma_j$:
$\Gamma_j\to\Lambda_j$ is a homeomorphism.
\enddemo

Set $J_0=\{j\in J:\ \Lambda_j\subset A\backslash Z,\ \dim \Lambda_j=d\}$
and $A'=\bigcup\{\Lambda_j:\  \Lambda_j\subset A,\ \dim \Lambda_j<d\}$. For
each
$j\in J_0$, $\Lambda_j$ is contained in some $T\in\cT$ of dimension
$d$. Thus $M_j$ is nonzero on $\Gamma_j$ and we have a
\L ojasiewicz inequality
$$
|M_j(\underline{b})|\ \ge\ C\cdot\left(\text{ dist}(\underline{b},\overline
\Gamma_j\backslash \Gamma_j)\right)^s \ ,
\leqno(6.2)
$$
for $\underline{b}\in\Gamma_j$, $j\in J_0$, where $s$ is a positive
integer.

Now take any integer $t\ge p+s$.
Suppose that $f\in\cC^t(\Omega)$ is $t$-flat on
$\vp^{-1}(A'\cup Z)$ and that, for each $a\in A$, there exists a
polynomial $W_a\in \BR[y]$ such that $T^tf(b,x)=W_a\left(\widetilde
T^t\vp(b,x)\right)\ \mod\  (x)^{t+1}$ for each $b\in\vp^{-1}(a)\cap E$.

Let $j\in J_0$.
For each $a\in\Lambda_j$, let
$V_j(a,y)\in\BR[y]$ be the unique polynomial (of degree $\le p$)
such that
$W_a(y)-V_j(a,y)\in\cR_p(a)$ and $\text{ supp}_y V_j(a,y)\subset
\Delta_j$.
Put
$$
V_j(a,y)=\sum_{\alpha\in\Delta_j} V_j^\alpha(a) y ^\alpha
\ ,
$$
$$
F_j(a,y)=\pi_l(V_j(a,y))=
\sum_{\alpha\in\Delta_j\atop|\alpha|\le l}V_j^\alpha(a) y^\alpha\ .
$$

\medskip

\proclaim{Lemma 6.3}
$\lim V_j^\alpha(a)=0$, when $a\to
\overline\Lambda_j\backslash\Lambda_j$.
\endproclaim

\demo{Proof}
Since $T^pf(b,x)=V_j\left(a,\widetilde T^p\vp(b,x)\right)\ \mod\ (x)^{p+1}$,
for
each
$b\in\vp^{-1}(a)\cap E$, we have
$$
D^\beta f(b_\nu)=\sum_{\alpha\in \Delta_j} V_j^\alpha
\left(\Phi(\underline{b})\right) L_\beta^\alpha (b_\nu)
$$
where $\underline{b}=(b_1,\ldots,b_q)\in\Gamma_j$, $\nu=1,\ldots,q$,
$\beta\in\BN^m$ and $|\beta|\le p$.
By Cramer's rule applied to this system,
$$
V_j^\alpha(\Phi(\underline{b}))\ =\ H_j^\alpha(\underline{b})\big/
M_j(\underline{b})
\leqno (6.4)
$$
for $\underline{b}\in\Gamma_j$, $\alpha\in\Delta_j$, where the $H_j^\alpha$
are $\cC^s$-functions on $\Omega^q$, $s$-flat on
$\overline\Gamma_j\backslash\Gamma_j\subset\Phi^{-1}\left(\overline\Lambda_j\backslash\Lambda_j\right)
\subset\Phi^{-1}(A')$. The lemma follows from (6.2).  $\square$
\enddemo
\medskip

\proclaim{Lemma 6.5}
$F_j$ is  a $\cC^l$-Whitney field on $\Lambda_j$.
\endproclaim

\demo{Proof} Let $\Psi_\nu$: $\Omega^q\to\BR^n$ and
$h_\nu$: $\Omega^q\to\BR$ ($\nu=1,\ldots,q$)
denote the mappings
$\Psi_\nu(\underline{b})=\Psi_\nu(b_1,\ldots,b_q)=\vp(b_\nu)$
and $h_\nu(\underline{b})=f(b_\nu)$.
Then $\Psi_\nu|E_\vp^q=\Phi$ and, for each $\underline{b}\in\Omega^q$,
the polynomial
$T^p\Psi_\nu(\underline{b},\underline{x})$ can be identified with
$T^p\vp(b_\nu,x)$.
Similarly, $T^ph_\nu(\underline{b},\underline{x})$ can be identified with
$T^pf(b_\nu,x)$.  Fix $a\in\Lambda_j$. Choose
$\underline{b}\in \Phi^{-1}(a)\cap E^q_\vp$ such that, for any
$W\in\BR[y]$, $W\in\cR_{p-1}(a)$ if and only if
$W\left(\widetilde T^p\vp(b_\nu,x)\right)=0\ \mod\ (x)^p$, $\nu=1,\ldots,q$.
Take
$S\in\cS^{(q)}$ such that $\underline{b}\in S$. Then $\Lambda_j$ is open in
$T=\Phi(S)$. We have
$$
T^p h_\nu(\underline{c},\underline{x}) \ =\
V_j\left(\Phi(\underline{c}),\widetilde
T^p\Psi_\nu(\underline{c},\underline{x})\right)
\ \mod\  (\underline{x})^{p+1}
\leqno(6.6)
$$
for all $\underline{c}\in S\cap\Phi^{-1}(\Lambda_j)$, $\nu=1,\ldots,q$.
By (6.4), the $V_j^\alpha$ are $\cC^1$ on $\Lambda_j$ (in fact,
$\cC^s$). Let $u$ be
any vector tangent to $\Lambda_j$ at $a$.
Since $\Phi|S$: $S\to T$ is a
submersion, there is $v\in T_{\underline{b}}S$ such that
$(d_{\underline{b}}(\Phi|S))(v)=\left(d_{\underline{b}}\Psi_\nu\right)(v)=u$.
By Lemma 3.1, Proposition 3.2 and (6.6),
$$
D_{a,u} V_j^{p-1}\left(a,\widetilde
T^p\Psi_\nu(\underline{b},\underline{x})\right) -
D_{y,u}V_j\left(a,\widetilde T^p\Psi_\nu(\underline{b},\underline{x})\right)
=0\ \mod\ (\underline{x})^p \ ;
$$
hence
$$
D_{a,u}V_j^{p-1}\left(a,\widetilde T^p\vp(b_\nu,x)\right) -
D_{y,u}V_j\left(a,\widetilde T^p\vp(b_\nu,x)\right)
=0\ \mod\ (x)^p \ .
$$
It follows that
$D_{a,u}V_j^{p-1}\left(a,y\right) -D_{y,u}V_j\left(a,y\right) \in\cR_{p-1}(a)$
and, by Lemma 5.4,
$$
D_{a,u}F_j^{l-1}\left(a,y\right) -
D_{y,u}F_j\left(a,y\right)
\in \pi_{l-1}(\cR_{p-1}(a))=\pi_{l-1}(\cR_p(a)).
$$
On the other hand,
$\text{ supp}_y\left[D_{a,u}F_j^{l-1}(a,y)-D_{y,u}F_j(a,y)\right]\subset
\Delta_j$;
thus $D_{a,u} F_j^{l-1}(a,y)
=D_{y,u}F_j(a,y)$. In virtue of
Proposition 3.2, this completes the proof of Lemma 6.5.  $\square$
\enddemo
\medskip

We can now finish the proof of Proposition 6.1.
Define $G\in\cC^0(A\cup Z)[y]$ by setting $G=F_j$ on
$\Lambda_j$, for each $j\in J_0$, and $G=0$ elsewhere.
Then $T^l f(b,x)=G(a,\widetilde T^l \vp (b,x))$ $\mod\  (x)^{l+1}$,
for each $a\in A\cup Z$ and $b\in\vp^{-1} (a)\cap E$.
By Proposition 3.4, the truncation $G^k(a,y)$ is a
$\cC^k$ Whitney field on $A\cup Z$.  $\square$


\bigskip

\Refs
\medskip


\ref\no 1
\by E. Bierstone, P.D. Milman,
\paper Composite differentiable functions,
\jour Ann. of Math. (2)
\vol 116
\yr (1982), \pages 541--558.
\endref

\ref\no 2
\bysame
\paper Algebras of composite differentiable functions,
\jour Proc. Symp. Pure Math.
\vol {\bf 40}
\yr (1983),
\pages 127--136.
\endref

\ref\no 3
\bysame
\paper  The Newton diagram of an analytic morphism, and
applications to differentiable functions,
\jour Bull. Amer. Math. Soc.
\vol {\bf 9}
\yr (1983),
\pages 315--318.
\endref

\ref\no 4 \bysame
\paper  Relations among analytic functions, I,
\jour Ann. Inst. Fourier (Grenoble)
\vol 37:1
\yr (1987),
\pages 187--239; and {\it II,}
{\bf 37:2}
(1987),
49--77.
\endref

\ref\no 5 \bysame
\paper  Semianalytic and subanalytic sets,
\jour Inst. Hautes \'Etudes Sci. Publ. Math.
\vol {\bf 67}
\yr (1988),
\pages 5--42.
\endref

\ref\no 6 \bysame
\paper {\it Geometric and differential properties of subanalytic sets,}
\jour Bull. Amer. Math. Sci.
\vol {\bf 25}
\yr (1991),
\pages 385--393.
\endref

\ref\no 7
\by E. Bierstone and G.W. Schwarz,
\paper Continuous linear division and extension of $\cC^\infty$-functions,
\jour Duke Math. J.
\vol {\bf 50}
\yr (1983),
\pages 233--271.
\endref

\ref\no 8
\by Z. Denkowska, S. \L ojasiewicz, J. Stasica,
\paper {\it Certaines properi\'et\'es \'el\'ementaires des ensembles
sous-analytiques,}
\jour Bull. Polish Acad. Sci. Math.
\vol {\bf 27}
\yr (1979),
\pages 529--536.
\endref

\ref\no 9
\by Z. Denkowska, S. \L ojasiewicz, J. Stasica,
\paper {\it Sur le th\'eor\`eme du compl\'ementaire pour les ensembles
sous-analytiques,}
\jour Bull. Polish Acad. Sci. Math.
\vol {\bf 27}
\yr (1979),
\pages 537--539.
\endref

\ref\no 10
\by G. Glaeser,
\paper {\it Fonctions compos\'ees diff\'erentiables,}
Ann. of Math. (2)
{\bf 77}
(1963),
193--209.
\endref

\ref\no 11
\by H. Grauert,
\paper {\it\"Uber die Deformation isolierter Singularit\"aten analytisher
Mengen,}
\jour Invent. Math.
\vol {\bf 15}
\yr (1972),
\pages 171--198.
\endref

\ref\no 12
\by R. H. Hardt,
\paper {\it Stratification of real analytic mappings and images,}
\jour Invent. Math.
\vol {\bf 28}
\yr (1975),
\pages 193--208.
\endref

\ref\no 13
\by R. H. Hardt,
\paper {\it Triangulation of subanalytic sets and proper light subanalytic
maps,}
\jour Invent. Math.
\vol {\bf 38}
\yr (1977),
\pages 207--217.
\endref

\ref\no 14
\by H. Hironaka,
\paper {\it Subanalytic sets},
\jour in: Number Theory, Algebraic Geometry and Commutative Algebra in
Honor of Y. Akizuki, Kinokuniya, Tokyo,
\yr 1973,
\pages 453--493.
\endref

\ref\no 15
\by H. Hironaka,
\paper {\it Introduction to the theory of infinitely near singular points,}
\jour Mem. Mat. Inst. Jorge Juan, No. 28,
Consejo Superior de Investigationes Cientificas,
Madrid,
\yr 1974.
\endref

\ref\no 16
\by S. \L ojasiewicz,
\book Ensembles semi--analytiques,
\publ Inst. Hautes \'Etudes Sci.,
\bookinfo Bures-sur-Yvette,
1964.
\endref

\ref\no 17
\by S. \L ojasiewicz,
\paper {\it Stratifications et triangulations sous-analytiques,}
\jour Univ. Stud. Bologna Sem. Geom.,
(1986),
\pages 83--89.
\endref

\ref\no 18
\by B. Malgrange,
\book {\it Ideals of Differentiable Functions,} Oxford Univ. Press, Bombay,
1966.
\endref

\ref\no 19
\by P.D. Milman,
\paper {\it The Malgrange-Mather division theorem,}
\jour Topology
\vol {\bf 16}
\yr (1977),
\pages 395--401.
\endref

\ref\no 20
\by W. Paw\l ucki,
\paper {\it On relations among analytic functions and geometry of subanalytic
sets,}
\jour Bull. Polish Acad. Sci. Math.
\vol {\bf 37}
\yr (1989),
\pages 117--125.
\endref

\ref\no 21
\by W. Paw\l ucki,
\paper {\it Points de Nash des ensembles sous-analytiques,}
\jour Mem. Amer. Math. Soc.
\vol{\bf 425}
\yr (1990).
\endref

\ref\no 22
\by W. Paw\l ucki,
\paper {\it On Gabrielov's regularity condition for analytic mappings,}
\jour Duke Math. J.
\vol {\bf 65}
\yr (1992),
\pages 299--311.
\endref

\ref\no 23
\by G.W. Schwarz,
\paper {\it Smooth functions invariant under the action of
a compact Lie group,}
\jour Topology
\vol {\bf 14}
\yr (1975),
\pages 63--68.
\endref

\ref\no 24
\by J. Cl. Tougeron,
\paper {\it An extension of Whitney's spectral theorem,}
\jour Inst. Hautes \'Etudes Sci. Publ. Math.
\vol {\bf 40}
\yr (1972),
\pages 139--148.
\endref

\ref\no 25
\by J. Cl. Tougeron,
\book Id\'eaux de Fonctions Diff\'erentiables,
\publ Springer,
\bookinfo Berlin-Heidelberg-New York,
1972.
\endref

\ref\no 26
\by J. Cl. Tougeron,
\paper {\it Fonctions compos\'ees diff\'erentiables: cas
alg\'ebrique,}
\jour Ann. Inst. Fourier (Grenoble)
\vol {\bf 30}
\yr (1980),
\pages 51--74.
\endref

\ref\no 27
\by H. Whitney,
\paper {\it Functions differentiable on the boundaries of regions,}
\jour Ann. of Math.
\vol {\bf 35}
\yr (1934),
\pages 482--485.
\endref

\ref\no 28
\by H. Whitney,
\paper {\it Differentiable even functions,}
\jour Duke Math. J.
\vol {\bf 10}
\yr (1943),
\pages 159--160.
\endref
\endRefs
\enddocument